\documentclass[12pt]{amsart}
\usepackage{graphicx}
\usepackage{authblk}
\usepackage{amsmath}
\usepackage{authblk}
\usepackage{amsthm}
\usepackage{amssymb}

\usepackage{tikz}
\usetikzlibrary{positioning,petri}
\usepackage{graphicx} 
\usepackage{pgfplots}
\usepackage{amsaddr}
\usepackage[margin=1in]{geometry}
\newcommand{\nonumberfootnote}[1]{%
  \begingroup%
  \def\thefootnote{}\footnote{#1}%
  \addtocounter{footnote}{-1}%
  \endgroup%
}

\usepackage{geometry}
\newtheorem{theorem}{Theorem}[section]
\newtheorem{corollary}{Corollary}[theorem]
\newtheorem{lemma}[theorem]{Lemma}

\newtheorem{Example}{Example}[theorem]
\newtheorem{proposition}[theorem]{Proposition}
\theoremstyle{definition}
\newtheorem{definition}{Definition}[section]

\usepackage{mathtools}
\DeclarePairedDelimiter{\ceil}{\lceil}{\rceil}
\theoremstyle{remark}

\begin{document}

\title{ Tropical Matrix Exponential
}
\author{ Ali M, Askar$^{a,\dagger}$ and  Mukherjee, Himadri$^a$}
\address{$^a$Department of Mathematics, BITS Pilani K. K. Birla Goa Campus, Goa, India}
\email{p20190037@goa.bits-pilani.ac.in (Askar$^\dagger$), himadrim@goa.bits-pilani.ac.in (Himadri) }
\nonumberfootnote{$\dagger$Corresponding author: Ali M, Askar}
\nonumberfootnote{AMS Classification 2010. Primary: 15A80 , 15A16}

\keywords{Max-plus algebra, Matrix exponential, Generalised eigenvector, Robust }

\maketitle

\noindent

\begin{abstract}
In this article, we introduce an exponential for tropical matrices and show that this series is essential for the analysis of certain kinds of stability in discrete event dynamic systems. A notion of a generalised eigenvector is introduced to discuss this kind of stability and prove it exists at most in the order of $1,p/2,p$, where $p$ is the period of the corresponding matrix. Thus characterizing the generalised eigenvectors of all powers of the matrix. Also, a sufficient condition is proved for the exponential of a matrix to be robust.  
\end{abstract}



\section{Introduction} 
Let $\mathbb{T}$ denote the tropical (max-plus) semi-field, defined as ($\mathbb{T} = \mathbb{R} \cup \{-\infty\}, \oplus, \otimes$), such that, for $a, b \in \mathbb{T},\  a \oplus b:= max \{a,b\}$ and $a \otimes b:= a+b$, the conventional sum in real line. As a convention, for the additional element $-\infty$, which will be denoted by $\varepsilon$ more often in this paper, we have, for any $a \in \mathbb{R},\  a \oplus \varepsilon = a$, and $a \otimes \varepsilon = \varepsilon$. \\
In many studies on max-plus algebra, significant focus has been directed towards the study
of power series of tropical matrices, motivated by their applications in the core problem of finding eigenvectors \cite{bapat-Gaubert,economics-note,max-at-work}. Among such power series, two series draw our attention, namely $\Gamma(A)$, the metric matrix and $\Delta(A)$, the Kleene star. For a given $A \in M_n(\mathbb{T})$, these are defined as follows. 
\[\Gamma(A):= A \oplus A^{(2)} \oplus A^{(3)} \oplus \cdots\]
and 
\[ \Delta(A):= I \oplus A \oplus A^{(2)} \oplus A^{(3)} \oplus \cdots \]
If the above infinite series $\Gamma(A)$ converges, then it gives the weights of heaviest paths of any length for all pairs of vertices, and hence it is called the metric matrix. It has been proved \cite{economics-note, thebook} that if the principal eigenvalue of $A$, $\lambda(A) \leq 0$, both $\Gamma(A)$ and $\Delta(A)$ are convergent, and in that case, they play a crucial role in finding eigenvectors of $A$, i.e., the solution to the eigenvalue-eigenvector problem,
\[A\otimes x = \lambda \otimes x, \text{ for }~ x \in \mathbb{T}^n- \{\varepsilon\},~ \lambda \in \mathbb{R}.\]
Moreover, when it's convergent, every column of $\Gamma(A)$ corresponding to the critical nodes (see Notations), are eigenvectors of $A$. \\

\begin{figure}[ht]
    \centering
 \resizebox{0.3\textwidth}{!}{
   \begin{tikzpicture}[yscale= -1.1, thick,>=stealth,
  every transition/.style={fill,minimum height=1mm,minimum width=3.5mm},
  every place/.style={draw,thick, scale= 0.7}]

        \node [transition, label=above left:$u^1$] (t1) {}; 
        \node [place, below =of t1, label=left:$g_1$] (p1) {}; 
        \node [transition, right=of t1, label=above:$u^2$] (t2) {}; 
        \node [place, below right =of t2, tokens=1, label=right:$g_2$] (p2) {}; 
         \node [transition, below=of p1, label=left:$x^1$] (t3) {}; 
        \node [place, right=of t3, tokens=1, label=below:$f_1$] (p3) {};
          \node [place, right=of p3, label=above:$f_1$] (p4) {};
         \node [transition, right=of p4, label=right:$x^2$] (t4) {}; 
        \node [place, below =of t3, label=left:$f_1$] (p5) {};
        \node [place, below =of p4, label=below:$f_2$] (p6) {};
        \node [place, below =of t4, tokens=1, label=below right:$h_2$] (p7) {};
        \node [transition, below =of p5, label=below:$x^3$ ] (t5) {};
        \node [place, left=of t5, tokens=2, label=left:$f_3$] (p8) {};
        \node [place, below=of p6,label= below left :$h_3$] (p9) {};
        \node [transition, below=of p9, label=left:$y$] (t6) {};
        
        \draw [->] (t1) -- (p1);
        \draw [->] (t2) -- (p2);
        \draw [->] (p1) -- (t3);
        \draw [->] (p2) -- (t4);
        \draw [->] (t3) -- (p5);
        \draw [->] (t4) to [bend right] (p3);
         \draw [->] (t4) -- (p6);
         \draw [->] (t5) -- (p9);
        \draw [->] (t5) to [bend left] (p8);
        \draw [->] (p8) to [bend left] (t5);
        \draw [->] (p5) -- (t5);
        \draw [->] (t3) to [bend right] (p4);
        \draw [->] (p4) to [bend right] (t4);

        \draw [->] (p6) to [bend left] (t5);
         \draw [->] (p3) to [bend right] (t3);
         \draw [->] (p9) -- (t6);
          \draw [->] (p7) to [bend right] (t6);
           \draw [->] (t4) -- (p7);
        
    \end{tikzpicture}
}
\label{Figure}
    \caption{Petri Net Example}
    \label{Fig.1}
\end{figure}
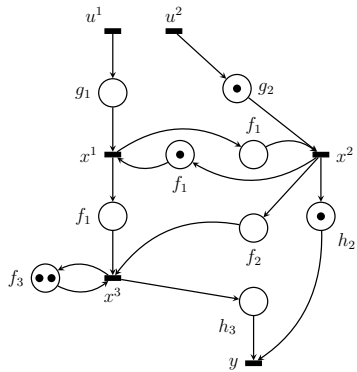

One of the major applications that Max-plus algebra found is in the discrete event dynamical systems (DEDS) (see for more details \cite{economics-note,thebook,max-at-work}). If $A$ is a realisation matrix for a DEDS, and $x(r)$ denotes the state variable at $r^{th}$ stage, then it follows a relation as follows:
\[x(r+1) = A \otimes x(r)\]
Then, the system reaches a steady regime if it progresses in regular steps after a stage. i.e., for some $\lambda \in \mathbb{R}$, and $r_0$,
\[x(r+1)= \lambda \otimes x(r) \text{ for all } r \geq r_0.\]
This is true, if and only if for some real number $\lambda$ and a positive integer $r$, $x(r)$ is a solution to the eigenvalue-eigenvector problem: $A \otimes x = \lambda \otimes x$. Hence, the eigenvalue-eigenvector problem is central to understanding these types of steady states. Stability or equilibrium analysis is one of the major aspects of studies in DEDS or in control theory. There are different types of stabilities that a system attains, and above discussed, one is one such stability. Readers may see \cite{khalil} for more notions of stabilities. In \cite{cohen} G. Cohen et al., and in \cite{stability}, Li Yanping et al. have defined stability in the following way:\\
The DEDS is stable if and only if there exists a
real number $\lambda$, and positive integers $k_0, d_0$ such that, 
\[x(k+d_0) = \lambda^{(d_0)}\otimes x(k), \text{ for all } k \geq k_0.\]
This happens,  if and only if for some real number $\lambda$ and positive integer $k$, $x(k)$ is a solution to the eigenvalue-eigenvector problem: 
\begin{equation}\label{gen.eigen}
    A^{(d_0)} \otimes x = \lambda^{(d_0)} \otimes x
\end{equation} 
\textit{Are they just eigenvectors of $A^{(d_0)}?$} In this article, we study solutions to \eqref{gen.eigen}, for all possible $d_0$'s, as a class of generalised eigenvectors, $GV(A)$ (including eigenvectors), which agrees with the definition given by G. Cohen et al. in \cite{cohen}. Giving credit to the periodic behaviour of max-plus matrices, the set $GV(A)$ encompasses generalized eigenvectors of the orders at most $1,~ p/2,~ p$, for some $p \in \mathbb{N}$.  This set is invariant for all powers and roots of $A$. i.e., $GV(A)= GV(A^{(k)})$, for every $k \in \mathbb{N}$. This will essentially classify eigenvectors of all powers and roots of $A$. Hence, it gives ample sense to study them together as generalised eigenvectors.\\
Petri nets are an important class of DEDS, where the above-mentioned notion of stability plays an important role. Consider the following example Figure \ref{Figure}, discussed in \cite{cohen}.

Where $f_i, g_i,h_i$  are times put on places, and $u^i ,x^i ,y$  are dates of transitions. Transitions $u^i$ are ``inputs" (tokens are brought from outside at dates $u_n^i$) and $y_n$ is the $n^{th}$ ``output" transition date. After the elimination of ``non-observable" variables \cite{petri-nets}, the system can be equivalently represented as follows.
\begin{align}
  \xi_n &= \xi_{n-1}F \oplus u_n G\\ 
  y_n &= \xi_nH
\end{align}
Then, G. Cohen et al. have proved in \cite{cohen} that, the transfer matrix $\mathcal{H}(z)= A_0 \oplus A_1 z^{-1} \oplus \cdots$, is realizable by a finite dimensional stable linear system of type $( 2 ) - ( 3 )$ if and only if there exist a $\lambda \in \mathbb{R},~ d_0, k_0 \in \mathbb{N}$, such that,
\[ \forall~ k \geq k_0,~~ A_{d_0+k} = \lambda^{(d_0)} \otimes A_{k}.\]  
When it is specialized for $A_{k}= A^{(k)}$, it goes back to \eqref{gen.eigen}.\\
In this article, we have developed a parallel approach for finding generalized eigenvectors of a matrix using the exponential of tropical matrices, $e^{(A)}$, in the same manner as $\Gamma(A)$  for finding eigenvectors of $A$ \cite{bapat-Gaubert}.
For a square matrix $A \in M_n(\mathbb{T})$, 
\[e^{(A)} := I \oplus (\frac{0}{1!} \otimes A) \oplus (\frac{0}{2!} \otimes A^{(2)}) \oplus \dots \oplus (\frac{0}{n!} \otimes A^{(n)}) \oplus \dots\]
Here, $n! := 1 \otimes 2\otimes \dots n = 1+2+ \dots +n = n(n+1)/2$, and $\frac{1}{n!}$ is the multiplicative inverse of $n!$, i.e., $-n(n+1)/2$.
Unlike the metric matrix $\Gamma $, the above series is proven convergent for every square matrix. The analysis of generalised eigenvectors using the exponential matrix holds significant utility across various stability analyses conducted on DEDS.\\ 

\textbf{Notations}\\
\begin{itemize}
    \item $I$ denotes the identity matrix of the appropriate size. $I = diag (0,0,\cdots,0)$, and every off-diagonal entries are $\varepsilon$.
    \item Powers: For $a,d \in \mathbb{R}$, $a^{(d)}$ denotes the Tropical power, and $a^d$ denotes the real(conventional algebra) power.\\
    Similarly, for a matrix $A\in M_n(\mathbb{T})$, $A^{(d)}$ denotes the Tropical power of $A$, and and $A^d$ denotes the real power.
    \item For $a,b \in ~ \mathbb{R}$, $\frac{a}{b} := a \otimes b^{(-1)}$, and $a/b$ denotes the conventional division.
    \item For a square matrix $A= [a_{ij}] \in M_n(\mathbb{T})$, the digraph associated with $A$ is denoted by $G_A$, with node set as $N= \{1,2, \cdots , n\}$, and there exist an edge from $i$ to $j$, with weight $a_{ij}$ if $a_{ij} > \varepsilon$.
    \item A square matrix $A$ is called irreducible if the corresponding digraph $G_A$ is strongly connected. i.e., There exists a directed path between any two vertices. 
    \item
If $\pi = (i_1,i_2),(i_2,i_3),\cdots,(i_k,i_1)$ is a cycle in $G_A$, then $(a_{i_1i_2}+a_{i_2i_3}+\cdots +a_{i_ki_1})/k$ is the corresponding cycle mean of $\pi$, and maximum cycle mean among all cycles of $G_A$ is denoted by $\mu(A)$. The cycles with a cycle mean equal to $\mu(A)$ are called critical cycles, and the union of all critical cycles is denoted by $C(A)$.
\item  Vertex set of $C(A)$ will be denoted by $E(A)$, and are called critical nodes.
\item For an irreducible square matrix $A$, $\lambda(A)$ denotes its unique eigenvalue, which is the maximum cycle mean in $G_A$\cite{perron-frobeniuos}. When it is clear from the context, $\lambda$ is used instead of $\lambda(A)$.
\end{itemize}
\section{Preliminaries}
Before we discuss the exponential series of tropical matrices, let us look into the power series of scalars from the tropical semi-field. For any $a, b\ \in\  \mathbb{T}$,\  $d_0: \mathbb{T} \rightarrow \mathbb{R}^+$, defined by $d_0(a,b) = |e^a - e^b|$, gives a metric in $ \mathbb{T} $. Then, it is easy to observe that, $\bigoplus_i a_i$, for $a_i \in \mathbb{T}$ converges, if and only if the corresponding sequence $\{a_i\}_i$ is bounded above. Moreover, if $\{a_i\}_i$ is bounded above, then it converges to the supremum.

\begin{definition}
    For any $a \in \mathbb{T}$, exponential map is defined as:
    \begin{equation}\label{scalar-exponential}
        e^{(a)}:= 0 \oplus \frac{a}{1!}  \oplus \frac{a^{(2)}}{2!}  \oplus \dots \oplus \frac{a^{(n)}}{n!}  \oplus \dots
    \end{equation}
\end{definition}

Note that, the exponential map has the following equivalent representation:
\begin{align}\label{redefine}
    e^{(a)} = \begin{cases}
		0 \text{  ; if } a \leq 1\\
                 \frac{1}{[a]!}\otimes a^{([a])}\text{ ;  if } a > 1	
		 \end{cases}
\end{align}
where $[\cdot]$ represents the greatest integer function. One can easily prove this fact using simple calculus techniques. As a consequence of this, we have the following result.
\begin{proposition}\label{factorial}
    For any positive integer $n \in \mathbb{N}$, $e^{(n )}= (n-1)!$. 
\end{proposition}
\begin{proof}
    From  \eqref{redefine}, we have,
     \begin{align*}
         e^{(n)} &= \frac{1}{n!}\otimes n^{(n)}\\
             &= \frac{1}{(n-1)!}\otimes n^{(n-1)}\\
             &= (n-1)!.
     \end{align*}
\end{proof}
 The above result is similar to the behaviour of $\Gamma$ function in conventional algebra. In fact, this is not surprising, as
 \begin{align*}
     e^{(a)} &= \bigoplus_{i=0}^\infty \frac{a^{(k)}}{k!}\\
             &= \int_{k \in \mathbb{N}\cup \{0\}}^\oplus  \frac{a^{(k)}}{e^{(k+1)}} dk \\
             &= \int_{k \in \mathbb{N}}^\oplus  (k-1)^a(e^k)^{-1} dk\\
             &= \Gamma(a).
 \end{align*} 
Readers may refer \cite{integral}, for a detailed understanding of integration and related analysis in tropical algebra. Thus, exponential for scalars gives an analogous function in tropical algebra for $\Gamma$ in conventional algebra. Readers may not confuse the $\Gamma$ series (metric matrix) with the $\Gamma$ function in conventional algebra over $\mathbb{C}$.
\begin{proposition}
   For $x,y \in \mathbb{T}$, if $x \oplus y=y$, then $e^{(x)} \oplus e^{(y)} = e^{(y)}$.
\end{proposition}
\begin{proof}
If $x \oplus y = y$, then we have $x \leq y$. Then, note that $x^{(m)} \leq y^{(m)},\ \forall m \ \in \mathbb{N}$. Now,
\begin{align*}
e^{(x)} &= 0 \oplus (\frac{1}{1!} \otimes x) \oplus (\frac{1}{2!} \otimes x^{(2)}) \oplus \dots \oplus (\frac{1}{n!} \otimes x^{(n)}) \oplus \dots  \\
&\leq 0 \oplus (\frac{1}{1!} \otimes y) \oplus (\frac{1}{2!} \otimes y^{(2)}) \oplus \dots \oplus (\frac{1}{n!} \otimes y^{(n)}) \oplus \dots\\
&= e^{(y)}.
\end{align*}
Hence, $e^{(x)} \oplus e^{(y)} = e^{(y)}$.
\end{proof}
We have the graph of $e^{(x)}$, which is piece-wise linear and mimics the graph of $e^x$, the exponential function in conventional algebra.
%
\[
\begin{tikzpicture}[scale=0.80]
\begin{axis}[
    axis lines = left,
    xlabel = \(x\),
    ylabel = {\(e^{(x)}\)},
]
\addplot [
    domain=0:1, 
    samples=50, 
    color=blue,
]
{0};
\addplot [
    domain=1:2, 
    samples=50, 
    color=red,
]
{x-1};
\addplot [
    domain=2:3, 
    samples=50, 
    color=blue,
    ]
    {2*x-3};
\addplot [
    domain=3:4, 
    samples=50, 
    color=red,
]
{3*x-6};
\addplot [
    domain=4:5, 
    samples=50, 
    color=blue,
]
{4*x-10};
\end{axis}
\end{tikzpicture}\]
The existence of the exponential map brings a natural question about the existence of the \textit{logarithmic function}, which inverts the exponential map $e^{(x)}$. From the above graph, it is evident that the exponential map is surjective to positive real numbers, and is one-one. This will enable us to define the $log$ function as follows.
\begin{definition}
    For $y>0$, $log (y)$ is defined as follows:
    \[ log(y) = (y\otimes1!)^{(1)} \oplus' (y\otimes2!)^{(1/2)} \oplus'  \cdots  \oplus' (y\otimes n!)^{(1/n)} \oplus' \cdots \]
where $\oplus'$ defined as the minimum, as for any $a, b \in \mathbb{T}$, $a \oplus' b := min \{a,b\}$.\\
\end{definition}
 For the above $log$, we can verify, for any $x> 0$, $log(e^{(x)})= x$. Moreover, the logarithm map has the following equivalent `continued fraction' form: 
\[log (y)~ =~ \frac{(y\otimes1!)^{(1)}\otimes \Big(\frac{(y\otimes2!)^{(1/2)} \otimes(\frac{(y\otimes3!)^{(1/3)} \otimes(\cdots) }{(y\otimes3!)^{(1/3)} \oplus (\cdots) }) }{(y\otimes2!)^{(1/2)} \oplus (\frac{(y\otimes 3!)^{(1/3)} \otimes(\cdots) }{(y\otimes3!)^{(1/3)} \oplus (\cdots) }) } \Big)}{(y\otimes1!)^{(1)}\oplus \Big(\frac{(y\otimes2!)^{(1/2)} \otimes(\frac{(y\otimes3!)^{(1/3)} \otimes(\cdots) }{(y\otimes3!)^{(1/3)} \oplus (\cdots) }) }{(y\otimes2!)^{(1/2)} \oplus (\frac{(y\otimes 3!)^{(1/3)} \otimes(\cdots) }{(y\otimes3!)^{(1/3)} \oplus (\cdots) }) } \Big)}\]   
    The proof is straightforward from the fact that, 
    \[a \oplus' b =  \frac{a \otimes b}{a \oplus b}\]
The graph of $log (x)$ is given below. Here as well, we can observe, that the tropical logarithmic map is piece-wise linear and mimics the logarithmic map in the conventional algebra. 
\[
\begin{tikzpicture}[scale=0.80]
\begin{axis}[
    axis lines = left,
    xlabel = \(x\),
    ylabel = {\(log(x)\)},
]
\addplot [
    domain=0:1, 
    samples=100, 
    color=blue,
]
{x+1};
\addplot [
    domain=1:3, 
    samples=100, 
    color=red,
]
{(x+3)/2};
\addplot [
    domain=3:6, 
    samples=100, 
    color=blue,
    ]
    {(x+6)/3};
\addplot [
    domain=6:10, 
    samples=100, 
    color=red,
]
{(x+10)/4};

\end{axis}
\end{tikzpicture}\]
\begin{definition}
Let $A \in M_{mn}(\mathbb{T})$, then $||A||_0 := e^{max \{a_{ij}\}}$, where maximum is taken over all entries of the matrix $A$, gives a norm on $M_{mn}(\mathbb{T})$. 
\end{definition}
\begin{lemma}
    Let $A \in M_{mp}(\mathbb{T})$ and $B \in M_{pn}(\mathbb{T})$, then $||A \otimes B||_0 \leq ||A||_0 ||B||_0$.
\end{lemma}
\begin{proof}
   We have, 
   \begin{align*}
        ||A \otimes B||_0 &= e^{max_{ij}\{\oplus_k a_{ik} +b_kj\}} \\
        &\leq e^{max\{a_{ij}\}+max\{b_{ij}\}} \\ 
        &= e^{max\{a_{ij}\}} e^{max\{b_{ij}\}}\\
        &= ||A||_0 ||B||_0.
   \end{align*}
\end{proof}
\begin{theorem}
Let $\{A_n\}_{n=1}^\infty$ be a sequence in $M_k(\mathbb{T})$. Then $\bigoplus_i A_i$ converges if and only if $\{A_n\}_{n=1}^\infty$ is bounded above.
\end{theorem}
\begin{proof}
    If $\{A_k\}_{k=1}^\infty$ is bounded above, then for each $i,j$, the sequence $(a_{n_{ij}})_{n=1}^\infty$ of $(i,j)^{th}$ entry of $A_n$, is bounded , hence has a supremum. Let $a_{ij}$ denotes the supremum of  $(a_{n_{ij}})_n$, and $A = (a_{ij})$ be the supremum matrix. Then, $\bigoplus_i A_i \ =\ A$. 

    Conversely, if $\bigoplus_i A_i$ is converging to some $A =  (a_{ij})$, then for any $n, i,j \in \mathbb{N}$, $a_{n_{ij}} \leq a_{ij}$. This gives,  $\{A_n\}_{n=1}^\infty$ is bounded above.
\end{proof}
\begin{lemma}\label{local_sum_of_products}
\begin{itemize}
    \item[a)]Let $a_i, b_i \in \mathbb{T}$, for $i = 1,2,\cdots n$. Then, $\bigoplus_{i=1}^n a_i b_i \leq \bigoplus_{i=1}^n a_i \bigoplus_{i=1}^n b_i $. 
    \item[b)]  Let $A_i, B_i \in M_n(\mathbb{T})$, for $i = 1,2,\cdots N$. Then, $\bigoplus_{i=1}^N A_i B_i \leq \bigoplus_{i=1}^N A_i \bigoplus_{i=1}^N B_i $.  
\end{itemize}
    
\end{lemma}
\begin{proof}
\begin{itemize}
    \item[a)]  Let $a_k := max\{a_i\}_{i=1}^n$ and $b_l := max\{b_i\}_{i=1}^n$. Then, for each $i$, we have the following. 
    \[ a_i  b_i \leq a_k  b_l\] 
    \[ \implies \bigoplus_{i=1}^n a_i  b_i \leq \bigoplus_{i=1}^n a_i  \bigoplus_{i=1}^n b_i \]
    \item[b)]  We have, 
    \begin{align*}
     \bigoplus_{i=1}^N A_i \bigoplus_{i=1}^N B_i &= \bigoplus_{i,j=1}^N A_i B_j\\
     & \geq \bigoplus_{i=1}^N A_i B_i
    \end{align*}
    
\end{itemize}
   
\end{proof}
 The following Corollary gives a comparison between any matrix power series and the well-established Gamma series. In particular, it gives a relation between the Gamma series and the exponential map. 
\begin{corollary}
    Let $A \in M_n(\mathbb{T})$, with $\lambda(A) > 0$, and $a_i \in \mathbb{T}$. Let $f(A):= \bigoplus_k a_k A^{(k)}$. Then, $f(A) \leq (\bigoplus_k a_k \lambda^{(k)}) \Gamma(\lambda^{-1}A)$. In particular, $e^{(A)} \leq e^{(\lambda)} \Gamma(\lambda^{-1}A)$.
 \end{corollary}
The periodic behaviour of powers of tropical matrices has been widely studied, notably in \cite{matrixpower, gen.period, period, power2}. P.Butkovic in \cite{matrixpower}, has given a characterisation for ``robust" matrices, with respect to the period of the matrix. In \cite{period}, M. Gavalec has given a formula and algorithm to find the period of an irreducible matrix. In this light, in the following section, we study the power series of matrices and its connection with the eigenvalue-eigenvector problem. Recall that, for $A \in M_n(\mathbb{T})$ and $x \in \mathbb{T}^n -\{\varepsilon\}$, the orbit of $A$ with starting vector $x$ is the sequence, $\mathcal{O}(A, x) = \{A^{(r)} \otimes x; r = 0, 1, \dots \}.$ Now,  $T(A) := \{x \in \mathbb{T}^n; \mathcal{O}(A, x) \cap V (A) \neq \varnothing\},$ where $V(A)$ is defined as collection of all eigenvectors of $A$. If $T(A) = \mathbb{T}^n - \{\varepsilon\}$, then $A$ is called robust. 
\begin{definition}
    We say that $A = (a_{ij} ) \in M_n(\mathbb{T})$ is ultimately periodic if there is a natural number $p$ such that the following holds for some $\lambda \in \mathbb{R}$ and $k_0$ natural:
$A^{(k+p)} = \lambda^{(p)} \otimes A^{(k)}$, for all $k \geq k_0$. If it is true for every $k \geq 1$, we say $A$ is strictly periodic.
If $p$ is the smallest natural number with this property, then we call $p$, the period of $A$ and denote it $per(A)$. If A is not ultimately periodic, then we set $per(A) = +\infty$.\\
Further, for this period $p$, the smallest natural number $k_0$, for which the relation, $A^{(k_0+p)} = \lambda^{(p)} \otimes A^{(k_0)}$ holds, is called the second order period, and is denoted by $Sper(A)$.
\end{definition}
\begin{theorem}\cite{matrixpower}\label{robust-charecter}
    Let $ A \in M_n(\mathbb{T})$ be irreducible. Then A is robust if and only if $per(A) = 1$.
\end{theorem}
\begin{theorem}\label{critical-cycle}\cite{period}
     Let $A = (a_{ij} ) \in M_n(\mathbb{T})$  be an irreducible matrix and $g_s$ be the $gcd$ of the lengths of critical cycles in the $s^{th}$ strongly connected component of $C(A)$. Then
$per(A) = lcm(g_1, g_2, . . .)$.
\end{theorem}

\section{Tropical Matrix Exponential }
\begin{definition}
    For a square matrix $A \in M_n(\mathbb{T})$, the exponential of $A$ is defined as, 
\[e^{(A)} := I \oplus (\frac{0}{1!} \otimes A) \oplus (\frac{0}{2!} \otimes A^{(2)}) \oplus \dots \oplus (\frac{0}{n!} \otimes A^{(n)}) \oplus \dots\]
Here, $n! := 1 \otimes 2\otimes \dots n = 1+2+ \dots +n = n(n+1)/2$, and $\frac{1}{n!}$ is the multiplicative inverse of $n!$, i.e., $-n(n+1)/2$.
\end{definition}
The following notations will be used throughout the paper.
\begin{itemize}
    \item For $A \in M_{m\times n}(\mathbb{T})$, $O(A) := \begin{cases}
    \ceil[\bigg]{  \bigoplus_{i = 1}^m \bigoplus_{j = 1}^n a_{ij}}, & \text{ if at least one } a_{ij}>0\\
    2, & \text{ if } a_{ij} \leq 0, \forall~ i,j 
     \end{cases}$\\
     where $\ceil{ a }$ is the least integer which is greater than $a$.
    \item For $A= (a_{ij})$, $a_{ij}^{k}$ denotes the $(ij)$ entry of $A^{(k)}$, unless stated otherwise.
\end{itemize}
\begin{lemma}
 Let $A \in M_n(\mathbb{T})$ be an irreducible matrix, then for any $k \in \mathbb{N}$,  
 $(\frac{0}{k!}\otimes A^{(k)} )\leq I \oplus \frac{0}{1!} \otimes A \oplus \cdots \oplus \frac{0}{t!} \otimes A^{(t)}$,
 where $ t = O(A)$.  
\end{lemma}
\begin{proof}
    If $\lambda(A) \leq 0$, then clearly $I \oplus \frac{0}{1!} \otimes A $ will dominate  $(\frac{0}{k!}\otimes A^{(k)} )$, for any $k \in \mathbb{N}$. Now, if $\lambda(A) >0$, then the entries in $\frac{0}{k!}\otimes A^{(k)}$ is bounded by $\frac{t^{(k)}}{k!} = k\cdot t - k \cdot (k+1)/2$, where $t = O(A)$. As $k$ increases large enough, this bound starts decreasing, and before decreasing, it attains a maximum at $k = t - \frac{1}{2}$. But as $k$ is an integer, let us take the maximum attained at the nearest integer to this, which is at $k = t$. Then we have, for any $k \in \mathbb{N},\  (\frac{1}{k!}\otimes A^{(k)} ) \leq I \oplus \frac{1}{1!} \otimes A \oplus \dots \oplus \frac{1}{t!} \otimes A^{(t)}$.  
\end{proof}
\begin{corollary}\label{speed-convergance}
For any matrix $A \in M_n(\mathbb{T})$, the infinite series expansion for $e^{(A)}$ converges. Further, the series will terminate in $O(A)$ steps. 
\end{corollary}
\begin{Example}
    The above discussion is interesting when the eigenvalue of the matrix $A$, $\lambda(A)\geq 0$. When $A$ is a negative matrix, i.e., $a_{ij} < 0$, for every $(i,j)$, the exponential power series terminates at the second step itself. The following example demonstrates it. \\
    Let $A = \begin{pmatrix}
        -2&-4&-1\\
        -3&-8&-4\\
        -1&-5&-6
    \end{pmatrix}$. Then, $e^{(A)} = \begin{pmatrix}
        0&-5&-2\\
        -4&0&-5\\
        -2&-6&0
    \end{pmatrix}$. For diagonal entries, entries from $I$ will dominate, and for every off-diagonal entry, entries from $\frac{A}{1}$ will dominate.\\
   \textbf{Remark:} Interestingly, for any matrix $A \leq 0$, i.e., $a_{ij} \leq 0, \forall ~ i,j$, the exponential $e^{(A)}$ is strongly definite. i.e., $\lambda(A)=0$, and $a_{ii}= 0$, for every $i$; and hence $e^{(A)}$ is robust. 
\end{Example}
\begin{Example}
    Let $B = \begin{pmatrix}
        3&1&2&4\\
        2&3&1&1\\
        4&2&2&1\\
        3&2&1&2
    \end{pmatrix}$. Then, $e^{(B)} = \begin{pmatrix}
        4&3&3&5\\
        3&3&1&3\\
        5&4&3&5\\
        4&3&2&4
    \end{pmatrix}$.
\end{Example}
    Any matrix $A \in M_n(\mathbb{T})$ can be transformed to its Frobenius normal form \cite{thebook, brualdi}, by simultaneous permutation of rows and columns. i.e., with appropriate permutation matrix $P$, $PAP^{-1} = \begin{pmatrix}
        A_{11}& \varepsilon & \cdots& \varepsilon\\
        A_{21}& A_{22}& \cdots & \varepsilon\\
        \cdots& \cdots& \cdots& \cdots\\
        A_{k1}&A_{k2} & \cdots & A_{kk}
    \end{pmatrix}$, where the diagonal blocks $A_{ii}$'s are irreducible square sub-matrices of $A$. 
\begin{lemma}
    Let $A \in M_n(\mathbb{T})$ be in its Frobenius normal form. i.e., $A = \begin{pmatrix}
        A_{11}& \varepsilon & \cdots& \varepsilon\\
        A_{21}& A_{22}& \cdots & \varepsilon\\
        \cdots& \cdots& \cdots& \cdots\\
        A_{k1}&A_{k2} & \cdots & A_{kk}
    \end{pmatrix}$. Then, the exponential of $A$ will have the form $ \begin{pmatrix}
        e^{(A_{11})}& \varepsilon & \cdots& \varepsilon\\
        *&e^{( A_{22})}& \cdots & \varepsilon\\
        \cdots& \cdots& \cdots& \cdots\\
        *& * & \cdots & e^{(A_{kk})}
    \end{pmatrix}$.
    \end{lemma}
    \begin{proof}
    The proof follows from the fact that $A^{(r)} = \begin{pmatrix}
        A_{11}^{(r)}& \varepsilon & \cdots& \varepsilon\\
        *& A_{22}^{(r)}& \cdots & \varepsilon\\
        \cdots& \cdots& \cdots& \cdots\\
        *& * & \cdots & A_{kk}^{(r)}
    \end{pmatrix} $.
    \end{proof}
    The above result enables us to check the properties of exponential of irreducible matrices, and those can be generalised to reducible matrices with the help of theories of reducible matrices. 
\begin{lemma}\label{eigenvalue}
    Let  Let $ A \in M_n(\mathbb{T})$ be irreducible, and $\lambda(A)$ be its eigenvalue, then $\lambda(e^{(A)})\ =\ e^{(\lambda(A))}$. 
\end{lemma}
\begin{proof}
    Let $ A \in M_n(\mathbb{T})$ be irreducible, then it has a unique eigenvalue \cite{perron-frobeniuos}, say $\lambda(A)$, with eigenvector $v (\neq \varepsilon)$. Then, we have $Av = \lambda v$.
\begin{align*}
    e^{(A)} \otimes v &= I \otimes v \oplus (\frac{0}{1!} \otimes A) \otimes v \oplus (\frac{0}{2!} \otimes A^{(2)})\otimes v \oplus \dots \oplus (\frac{0}{n!} \otimes A^{(n)}) \otimes v \oplus \dots\\
    &= v \oplus \frac{0}{1!} \otimes \lambda \otimes v \oplus \frac{0}{2!} \otimes \lambda^{(2)} \otimes v \oplus \dots \oplus \frac{0}{n!} \otimes \lambda^{(n)} \otimes v \oplus \dots\\
    &= e^{(\lambda)} \otimes v 
\end{align*}
This gives $v$ is an eigenvector for $e^{(A)}$, with respect to the eigenvalue $e^{(\lambda)}$.  
\end{proof}
\begin{lemma}
    For $X, Y \in M_n(\mathbb{T})$, if $X \oplus Y = Y$, then $e^{(X)} \oplus e^{(Y)} = e^{(Y)}$.
\end{lemma}
\begin{proof}
Let $X = (x_{ij}),\ Y = (y_{ij}) \in M_n(\mathbb{T})$.
    If $X \oplus Y = Y$, then for every $i,j$, $x_{ij} \leq y_{ij}$. Let $x_{ij}^{k}$ denotes the $ij^{th}$ entry of $X^{(k)}$. Then, note that, $x_{ij}^{k}$ is the maximum weight of a walk on $G_X$, the digraph induced by $X$, connecting $i$ to $j$ with length $k$. Since  $x_{ij} \leq y_{ij},\ \forall i,j$, we have, on each walk, the maximum weight of $X$ is less than or equal to that of $Y$. Hence, $X^k \oplus Y^k = Y^k$. This gives, 
    \begin{align*}
        e^{(X)} &= I \oplus (\frac{0}{1!} \otimes X) \oplus (\frac{0}{2!} \otimes X^{(2)}) \oplus \dots \oplus (\frac{0}{n!} \otimes X^{(n)}) \oplus \dots\ \\
         &\leq I \oplus (\frac{0}{1!} \otimes Y) \oplus (\frac{0}{2!} \otimes Y^{(2)}) \oplus \dots \oplus (\frac{0}{n!} \otimes Y^{(n)}) \oplus \dots \\ 
         &= e^{(Y)}
    \end{align*}
This gives, $e^{(X)} \oplus e^{(Y)} = e^{(Y)}$.
\end{proof}
\begin{corollary}
For a matrix $A \in M_n(\mathbb{T})$, let $H_A$ denotes the cone $\{X\ \in M_n(\mathbb{T}) | X \oplus A = A \}$. Then, $e^{(H_A)}:= \{e^{(X)} | X \in H_A \}\ \subseteq H_{e^{(A)}}$. 
\end{corollary}
\begin{proposition}
    Let $A , B\in M_n(\mathbb{T})$, such that $A = \mu B$, for some scalar $\mu>1$, and $\lambda(B)=0$. Further, assume $A$ is strictly periodic, robust matrix. Then, $e^{(A)} = e^{(\mu)} B$.
\end{proposition}
\begin{proof}
   First note that, $\lambda(A)= \mu$. Now,
    \begin{equation*}
        \begin{aligned}
            e^{(A)} &= I \oplus \frac{A}{1!}\oplus \frac{A^{(2)}}{2!}\oplus \cdots  \\
                   &= \frac{A}{1!}\oplus \frac{A^{(2)}}{2!}\oplus \cdots, \text{ since }\mu>1 \\
                   &= \frac{A}{1!}\oplus \frac{\mu A}{2!}\oplus \frac{\mu^{(2)} A}{3!} \cdots\\
         &= \frac{\mu B}{1!}\oplus \frac{\mu^{(2)} B}{2!}\oplus \frac{\mu^{(3)} B}{3!}\cdots \\
         &= B (\frac{\mu }{1!}\oplus \frac{\mu^{(2)} }{2!}\oplus \frac{\mu^{(3)} }{3!}\cdots)\\
         &= e^{(\mu)} B
        \end{aligned}
    \end{equation*}
\end{proof}
For the above-mentioned class of matrices, i.e., strictly periodic robust matrices, exponential preserves the structure of the original matrix. If we see $B$ as a skeleton of $A$, we have got $e^{(A)}$ has the same skeleton. 
\begin{proposition}
\begin{itemize}
    \item[a)] Let $a > 0$. Then, $e^{(1\otimes e^{(a)})}= e^{(a)} e^{(e^{(a)})}$.
    \item[b)] Let $A \in M_n(\mathbb{T})$, such that $\lambda(A)>1$. Then, $e^{(1\otimes e^{(A)})}= e^{(A)} e^{(e^{(A)})} $.
\end{itemize}    
\end{proposition}
\begin{proof}
\begin{itemize}
    \item[a)] Let $n \leq e^{(a)} \leq n+1$. Then by Lemma \ref{redefine}  we have, $e^{(e^{(a)})}= \frac{n e^{(a)}}{n!}$. This gives,
    \begin{equation*}
        \begin{aligned}
            e^{(1\otimes e^{(a)})} &= \frac{(n+1) e^{(a)\otimes 1}}{(n+1)!}\\
            &= \frac{n e^{(a)}}{n!} \otimes e^{(a)} =  e^{(e^{(a)})} e^{(a)}.
        \end{aligned}
    \end{equation*}
    \item[b)]  Let $e^{(A)} = B$. Then,
    \begin{equation*}
        \begin{aligned}
            e^{(1\otimes e^{(A)})} &= e^{(1\otimes B)}\\
             &= \frac{1\otimes B}{1!}\oplus \frac{2\otimes B^{(2)}}{2!}\oplus \cdots\\
             &= B (I \oplus \frac{B}{1!}\oplus \frac{B^{(2)}}{2!}\oplus \cdots)\\
             &= B e^{(B)}\\
             &= e^{(A)} e^{(e^{(A)})} 
        \end{aligned}
    \end{equation*}
\end{itemize} 
\end{proof} 
\subsection{Robust Exponentials}
In the following results, we will examine when the exponential of a matrix is robust. The study of robust matrices is particularly important in discrete-event dynamic systems \cite{discrete-event-1,discrete-event-2,gaubert}. P. Butkovic and R.A. Cuninghame-Green in \cite{matrixpower} have proved that the sequence of Markov parameters of the discrete-event dynamic system $(A, b, c)$ is
ultimately linear if and only if $A$ is robust.  

\begin{lemma}
     Let $A \in M_n(\mathbb{T})$ be an irreducible matrix, such that $A$ is periodic with period 1. Then columns of $A$, $a_j \in T(e^{(A)}), \forall j$. 
\end{lemma}
\begin{proof}
Let A be periodic with period 1. Then, there exist a $k_0$, natural number, such that $A^{(k+1)} = \lambda A^{(k)},\ \forall k \geq k_0$. We have, 
 \begin{align*}
 e^{(A)} &=\ I \oplus  (\frac{1}{1!} \otimes A) \oplus (\frac{1}{2!} \otimes A^{(2)} \oplus \dots \oplus (\frac{1}{n!} \otimes A^{(n)}) \oplus \dots \\
 &= I \oplus  (\frac{1}{1!} \otimes A) \oplus (\frac{1}{2!} \otimes A^{(2)} \oplus \dots \oplus (\frac{1}{t!} \otimes A^{(t)} )\\
 \implies e^{(A)} \otimes A^{(k)} &= A^{(k)} \oplus  (\frac{1}{1!} \otimes A^{(k+1)}) \oplus (\frac{1}{2!} \otimes A^{(k+2)} \oplus \dots \oplus (\frac{1}{t!} \otimes A^{(k+t)})\\
 &= (1 \oplus  \frac{\lambda}{1!} \oplus \frac{\lambda^{(2)}}{2!}  \oplus \dots \oplus \frac{\lambda^{(t)}}{t!}) \otimes A^{(k)}\\
 &= e^{(\lambda)} \otimes A^{(k)}.
 \end{align*}
\end{proof}
\textbf{Note:} Let $i,j \in E(A)$, and we say $i,j$ belongs to same equivalent class of vertices if they belong to the same critical cycle. In \cite{matrixpower}, P.Butkovic et al. have proved that for any strongly irreducible matrix (i.e., all powers of $A$ are irreducible), $E(A) = E(A^{(k)})$, and equivalent classes of $E(A^{(k)})$ are either same as equivalent classes of $E(A)$ or are their refinements. Also, the only refinements happening in powers are singleton sets, corresponding to the self-loops.  
\begin{theorem}
    Let $ A \in M_n(\mathbb{T})$ be strongly irreducible and robust. Then $e^{(A)}$ is robust.
\end{theorem}
\begin{proof}
    In \cite{commuting-same-eigenspace}[Theorem 5.5], R.D. Katz et al. have proved that two commuting irreducible matrices have the same critical nodes. i.e., for irreducible matrices $A,B$, if $AB=BA$, then $E(A)= E(B)$. Note that, $A$ and $e^{(A)}$ are commuting matrices. Hence, we have $E(A)= E(e^{(A)})$. More over, $C(e^{(A)})$ is a refinement of $C(A)$. Now, if $A$ is robust, then by Theorem \ref{critical-cycle}, $g_s$ is 1 in every strongly connected component, where $g_s$ is the $gcd$ of the lengths of critical cycles in the $s^{th}$ strongly connected component of $C(A)$. On refinement, this $gcd$ doesn't change. Hence, period of $e^{(A)}$ is 1, and $e^{(A)}$ is robust.  
\end{proof}
\begin{lemma}
    Let $ A \in M_n(\mathbb{T})$ be irreducible. If there exists an $i_s$ corresponding to every strongly connected component of $C(A)$,  such that, the $(i_si_s)^{th}$ entry of A, $(a_{i_si_s}) = \lambda(A)$, then $A$ is robust. 
\end{lemma}
\begin{proof}\label{diagonal-entry}
    Without loss of generality, let us assume $C(A)$ is strongly connected, as we can examine the hypothesis for each strongly connected component. Then, if there exist an $i_s$, such that $(a_{i_si_s}) = \lambda(A)$, then the self-loop around $i_s$, of length 1, is a critical cycle. We have, $gcd (\text{ lengths of all critical cycles of } C(A)) = 1$. Then by Theorem \ref{critical-cycle}, we have $A$ is robust.
\end{proof}
\begin{lemma}\label{accumulation}
   Let $ A \in M_n(\mathbb{T})$ be irreducible. If there exists a critical cycle of length $m$, $\sigma = (i_0, i_1, \dots, i_{m-1}, i_0 )$, then for any $k \in \mathbb{N}$, such that $m$ divides $k$, $(a_{i_0i_0}^{k})= \lambda(A)^{(k)}$. Further, if $C(A^{(k)})$ is strongly connected, then $A^{(k)}$ is robust.
\end{lemma}
\begin{proof}
    Let $i_0 \in E(A)$, such that there exists a critical cycle $\sigma = (i_0, i_1, \dots, i_{m-1}, i_0 )$ of length $m$. Take $k \in \mathbb{N}$ such that $m$ divides $k$, and say $k/m = d$. Then, $(i_0i_0)$ entry of $A^{(k)},\ (a_{i_0i_0}^{(k)})$, is the maximum weight of a walk of length $k$, starting from node $i_0$ and ending at $i_0$. The walk by repeating $\sigma$, $d$ times gives a walk of length $k$, and its weight is $d \cdot\ weight(\sigma)$. But, as $\sigma$ is a critical cycle of length $m$, we have  $weight(\sigma)= m \cdot \lambda(A)$. Thus, the above-mentioned walk from $i_0$ to $i_0$ should have a weight $d \cdot m \cdot \lambda(A) = \lambda(A)^{(k)} $. Since $\lambda(A)^{(k)} $  is the maximum possible weight in $A^{(k)}$, we get $(a_{i_0i_0}^{k}) = \lambda(A)^{(k)} $. Also, we have $\lambda(A^{(k)})= \lambda(A)^{(k)}$. Now, if $C(A^{(k)})$ is strongly connected, we have $A^{(k)}$ is robust, by Lemma \ref{diagonal-entry}.     
\end{proof}
\begin{lemma}
    Let $ A \in M_n(\mathbb{T})$ be irreducible, and $i_0 \in E(A)$. Let $m_1,m_2, \dots, m_p$ be the lengths of critical cycles containing $i_0$. Then for $k \in \mathbb{N}$, such that for any real-linear combinations $m_j$s do not divides k, we have $(a_{i_0i_0}^{k}) < \lambda(A)^{(k)}$.
\end{lemma}
\begin{proof}
    First, note that $(a_{i_0i_0}^{k}) \leq \lambda(A)^{(k)}$, otherwise, the self-loop of length 1 around $i_0$ has weight greater than $\lambda(A)^{(k)}$, which contradicts the fact that $\lambda(A^{(k)})= \lambda(A)^{(k)}$. Now, if $(a_{i_0i_0}^{k}) = \lambda(A)^{(k)}$, then there exists a combination of critical cycles which are beginning and ending at $i_0$ and have a total length $k$. But this contradicts our assumption that $k$ is not a multiple of any real-linear combination of $m_j$s. Thus, we have  $(a_{i_0i_0}^{k}) < \lambda(A)^{(k)}$.
\end{proof}
\begin{theorem}
     Let $ A \in M_n(\mathbb{T})$ be irreducible, then we have the following. 
     \begin{itemize}
         \item[a)]  If every strongly connected component of $C(A)$ has a critical cycle such that its length divides $[\lambda(A)]$, then $e^{(A)}$ is robust. 
         \item[b)]  If $\lambda(A)$ is a positive integer, and if every strongly connected component of $C(A)$ has a critical cycle such that its length divides $\lambda(A)$ or $\lambda(A)-1$, then $e^{(A)}$ is robust.
     \end{itemize}
\end{theorem}
\begin{proof}
    Without loss of generality, let us assume $C(A)$ is strongly connected. Take $[\lambda(A)] = r$. Now, if there exist an $i_0 \in E(A)$, such that length of a critical cycle around $i_0$ divides $r$, then by Lemma \ref{accumulation}, we have $a_{i_0i_0}^{r} = \lambda(A)^{(r)}$. But we have, $\lambda(e^{(A)}) = e^{(\lambda(A))}$, by Lemma \ref{eigenvalue}. Now, using Lemma \ref{redefine}, we get $e^{(\lambda(A))} = \frac{1}{r!} \otimes \lambda(A)^{(r)}$. Hence the maximum value that can be attained by any diagonal entry of $e^{(A)}$ is  $\frac{1}{r!} \otimes \lambda(A)^{(r)}$. But, in the expression, 
    \[e^{(A)} := I \oplus (\frac{1}{1!} \otimes A) \oplus (\frac{1}{2!} \otimes A^{(2)}) \oplus \dots \oplus (\frac{1}{r!} \otimes A^{(r)}) \oplus \dots\]
    $(i_0i_0)$ entry is attaining this value at the $r^{th}$ term, and hence $(e^{(A)})_{i_0i_0} = \frac{1}{r!} \otimes \lambda(A)^{(r)}$. Then by Lemma \ref{diagonal-entry}, $e^{(A)}$ is robust.\\
 Further, if $\lambda(A)$ is a positive integer, then, by the Lemma \ref{factorial}, \[e^{(\lambda(A))} = \frac{1}{\lambda(A)!} \otimes \lambda(A)^{(\lambda(A))} = \frac{1}{(\lambda(A)-1)!} \otimes \lambda(A)^{(\lambda(A)-1)}.\]
Hence, if the length of any critical cycle divides $\lambda(A)$ or $\lambda(A)-1$, then $e^{(A)}$ is robust.
\end{proof}
The above result characterizes when the exponential of a matrix is robust. Here, even though $A$ is not robust, $e^{(A)}$ is robust, with the same eigenvectors of $A$. Hence, we are getting a robust realization matrix with the same eigenvectors (state vector corresponding to an equilibrium state) as the previous system.
\begin{Example}
    Let $A= \begin{pmatrix}
        4&3&2\\
        5&2&6\\
        3&4&2
    \end{pmatrix}$. Then, the corresponding digraph $G_A$, and the critical digraph $C(A)$ are as follows respectively:

    \[ \includegraphics[width=6cm, height=5cm]{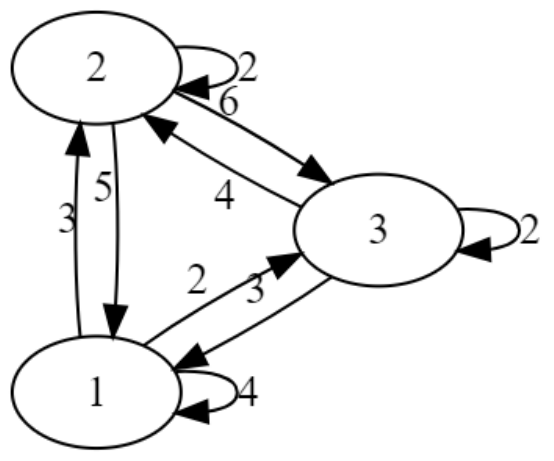} \includegraphics[width=2.5cm, height=4cm]{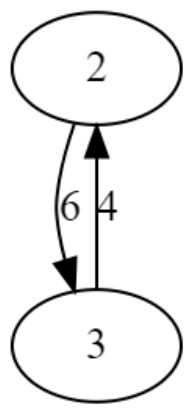} \]
    For this given $A$, the exponential $e^{(A)}= \begin{pmatrix}
        8&8&9\\
        10&10&11\\
        9&9&10
    \end{pmatrix}$. Here, we can see, $\lambda(A) = 5$, and length of critical cycle is $2$, and which divides $4 = \lambda(A)-1 $. One can verify $e^{(A)}$ is robust, using Theorem \ref{robust-charecter}  and Theorem \ref{critical-cycle}.  
\end{Example}

\section{Generalised Eigenvectors}
\begin{definition}
     Let $A \in M_n(\mathbb{T})$ and $x \in \mathbb{T}^n -\{\varepsilon\}$. Then we say $x$ is a generalised eigenvector of $A$, of order $m \in \mathbb{N}$ and corresponding to the eigenvalue $\lambda$, if,
     \[A^{(m)} \otimes x = \lambda^{(m)} \otimes x, \text{ but, } A^{(m-1)} \otimes x \neq \lambda^{(m-1)} \otimes x. \]
    Set of all such generalised eigenvectors is denoted by $GV(A) := \{x  \in \mathbb{T}^n -\{\varepsilon\}|\  A^{(m)} \otimes x = \lambda^{(m)} 
\otimes x, m=1,2, \dots \}.$
\end{definition}
\begin{definition}
   Let $A \in M_n(\mathbb{T})$. Define $GT(A) := \{x \in \mathbb{T}^n;\  O(A, x) \cap GV (A) \neq \varnothing\}$. We say $A$ is quasi-robust, if $GT(A) = \mathbb{T}^n - \{\varepsilon\}$. 
\end{definition}
\begin{theorem}
    Let $A \in M_n(\mathbb{T})$ be irreducible, then $A$ is quasi-robust if and only if $A$ is periodic, with $per(A)= p < \infty$.   
\end{theorem}
\begin{proof}
  Let A be a periodic matrix with $per(A) = p$. Then, there exist $k_0 \in \mathbb{N}$ and $\lambda \in \mathbb{R}$, such that for all $k \geq k_0$, $A^{(k+p)} = \lambda^{(p)} A^{(k)}$. Now, if $x \in \mathbb{T}^n -\{\varepsilon\}$, then, $A^{(k+p)} \otimes x\  =\ \lambda^{(p)}\  A^{(k)}\  \otimes x.$ This gives $A^{(k)} \otimes x \in GV(A)$, thus A is quasi-robust. 

Now, if $A$ is quasi-robust, then for $x = c_j$, the $j^{th}$ column of A, there exist a $k_j, m_j \in \mathbb{N}$, such that for all $k \geq k_j$, $A^{k+m_j} \otimes c_j = \lambda^{(m_j)} A^{(k)} c_j$, for some $\lambda$. Take $k_0 = max \{k_1,k_2,\dots, k_n\}$, and $m_0 = lcm(m_1,m_2, \dots, m_n)$. Then for any $k \geq k_0$, we have $A^{(k+m_0)}= \lambda^{(m_0)} \otimes A^{(k)}$. This implies, $A$ is periodic with $per(A) = m_0$. 
\end{proof}
\begin{proposition}
    Let $A,B$ commute, then they have a common generalised eigenvector. 
\end{proposition}
\begin{proof}
    In \cite{commuting-same-eigenspace}, R. D. Katz et al. have proved that if $A, B$ commute, then they have a common eigenvector. The rest is followed by the fact that, if $A, B$ commute, then $A^{(k)}, B^{(m)}$ commute for any $m,k \in \mathbb{N}$.
\end{proof}
\begin{lemma}\label{p-gen.eigenvector}
     Let $ A \in M_n(\mathbb{T})$ be irreducible, with $per(A) = p >1$. Then, $A$ has a generalised eigenvector of order $p$
\end{lemma}
\begin{proof}
    Since $A$ has a period, which is greater than 1, $A$ is not robust. i.e., there exists a non-trivial vector $v$, such that $v \notin T(A)$. Also, we have, for some some $k_0 \in \mathbb{N}$, and for every $k \geq k_0$ , $A^{(p+k)} = \lambda^{(p)} A^{(k)}$,as $A$ is periodic. Let $k_0$ be the least integer for which this happens. Then, for the above $v$, we have,
    \[A^{(p+k_0)}v = \lambda^{(p)} A^{(k_0)}v\]
    This implies, $A^{(k_0)}v$ is a generalised eigenvector of order $p$ for $A$. 
\end{proof}
\begin{theorem}
    Let $ A \in M_n(\mathbb{T})$ be irreducible, with $per(A) = p$. The only orders of generalized eigenvectors of $A$ are $p$ and $p/2$, if $p$ is even, and it is $p$ alone if $p$ is odd. 
\end{theorem}
\begin{proof}
Let $GT_{m_i}$ be the collection of all vectors $x$, which eventually satisfy the equation $A^{(m_i)}x =\lambda^{(m_i)} x$. i.e., There exists a $k_0 \in \mathbb{N}$, such that for every $k \geq k_0$ , $A^{(m_i+k)}x = \lambda^{(m_i)} A^{(k)}x$. Note that, $GT_{m_i} \cap GT_{m_j} \neq \varnothing$, if and only if $m_i | m_j$.\\
   Now, let $x_0 \in \mathbb{T}^n - \{\varepsilon\}$ be a generalised eigenvector of $A$ of order $m$, corresponding to the eigenvalue $\lambda$. Then we have, $A^{(m)}x_0 =\lambda^{(m)} x_0 $. Since, $A$ is periodic with period $p$, we have, for some $k_0 \in \mathbb{N}$, for every $k \geq k_0$,
    \[A^{(p+k)}x_0 = \lambda^{(p)} A^{(k)}x_0\]
    This gives $x_0 \in GT_p$. But since $x_0$ is a generalized eigenvector of order $m$, we have $x_0 \in GT_m$. This implies if $m \leq p$, $m$ should be a divisor of $p$.  If $m > p$, then $m$ must be a multiple of $p$. But, since $x_0$ is already cyclic with the order $p$, it contradicts the fact that the order of generalized eigenvector is the least integer with respect to which the vector is cyclic. Hence, $m$ can not exceed $p$.\\
    Now, Let $v \in GT_{m_1}\cap GT_{m_2}$, for some $m_1 < m_2$. Then there exists, $k_1,k_2$, such that, 
    \[A^{(k_1+m_1)} v = \lambda^{m_1}A^{(k_1)}v \text{ and } A^{(k_2+m_2)} v = \lambda^{m_2}A^{(k_2)}v \]
Then, we have, 
\begin{equation*}
    \begin{aligned}
        \lambda^{(m_1)}A^{k_1+(m_2-m_1)}v &=  A^{(m_2-m_1)}A^{(k_1+m_1)}v\\
        &= A^{(k_1+m_2) }v\\
        &= A^{(k_2+m_2)}A^{(k_1-k_2)} v \\
        &= \lambda^{m_2}A^{(k_1)}v
    \end{aligned}
\end{equation*}
Let $m':= m_2-m_1$. Then, from the above equations, it follows that,
\[A^{(k_1+m')}v = \lambda^{m'}A^{k_1}v.\]
This gives, $v \in GT_{m'}$. i.e., $v$ is cyclic with respect to the order $m'$, which is smaller than $m_2$.\\
By Lemma \ref{p-gen.eigenvector}, $p$ is already an order for a generalized eigenvector. Now, if $m$ is an order for a generalised eigenvector, then $m$ must divide $p$. This gives, $GT_p \cap GT_m \neq 0$. If $v \in GT_p \cap GT_m$, then by the above claim,  $v \in GT_{p-m}$. Then, this gives either $v$ is an eigenvector or $p-m$ divides $p$. Since $m$ is an order of generalized eigenvector, we can find some vector $v$ in $GT_p \cap GT_m$, which is not an eigenvector. Hence, $p-m$ should divide $p$. But, both $m$ and $p-m$ divides $p$, only if $p$ is even and $m = p/2$.   
\end{proof}
\textbf{Remark:} The above result classifies all eigenvectors of powers of $A$, and roots of $A$ as well to that matter. They are not very different classes. Hence, studying them together will be more useful, and the term `generalized eigenvectors' makes more sense in this aspect. 
\begin{Example}
    Consider the following matrix, $A = \begin{pmatrix}
        2&0&-1&3&1\\
        3&-1&1&2&0\\
        0&4&-1&2&1\\
        1&2&2&1&0\\
        -1&0&1&0&0
    \end{pmatrix}$
Then, $\lambda(A)= 3$. With respect to this eigenvalue,  $v = \begin{pmatrix}
    0\\
    -1\\
    1\\
    -1\\
    -2
\end{pmatrix} $ is a generalised eigenvector of order 2, and for $A$, and $u = \begin{pmatrix}
    0\\
    -1\\
    0\\
    -1\\
    -2
\end{pmatrix}$ is a generalised eigenvector of order 4.
\end{Example}
\begin{theorem}\label{gen.eigenvector relation}
    Let $ A \in M_n(\mathbb{T})$ be irreducible, with $per(A) = p$, and second order period $Sper(A)= N_0$. Further, assume that $a_{ii} \geq 1$. Then we have the following:
    \begin{itemize}
        \item[(a)] Let $A$ be strictly periodic. Then, every column of $e^{(A)}$ is a generalised eigenvector for $A$. 
        \item[(b)] If $N_0+2p < O(A)$, then every columns of $e^{(A)}$, are generalised eigenvectors of $A$.  
    \end{itemize}   
\end{theorem}
\begin{proof}
    \begin{itemize}
        \item[(a)] Since $a_{ii} \geq 1$, exponential can be written as 
        \[e^{(A)} = \frac{A}{1!}\oplus \frac{A^{(2)}}{2!}\oplus \cdots \oplus \frac{A^{(t)}}{t!}, \text{  where } t= O(A). \]
    Then, \begin{equation*}
        \begin{aligned}
            A^{(p)}e^{(A)} &= \frac{A^{(p+1)}}{1!}\oplus \frac{A^{(p+2)}}{2!}\oplus \cdots \oplus \frac{A^{(p+t)}}{t!}\\
            &= \lambda^{(p)} \Big(\frac{A}{1!}\oplus \frac{A^{(2)}}{2!}\oplus \cdots \oplus \frac{A^{(t)}}{t!} \Big)\\
            &= \lambda^{(p)} e^{(A)}
        \end{aligned}
    \end{equation*}
    \item[(b)] First, note that we can change $O(A)$ by multiplying with a suitable scalar, and it does not change the $per(A)$ or $Sper(A)$. Hence, this result can be applied to any matrix after rescaling appropriately. Also, every columns of $A^{(N_0+p+k)}, \text{ for } k \in \mathbb{N}$, are generalised eigenvectors of $A$, as $ A^{(N_0+p+k)} = \lambda^{(p)} A^{(N_0+k)}$.\\
    Now, we know that $e^{(A)}$ takes at most $O(A) = t$ steps to converge. i.e.,
    \begin{equation*}
        \begin{aligned}
           e^{(A)} &= \frac{A}{1!}\oplus \frac{A^{(2)}}{2!}\oplus \cdots \oplus \frac{A^{(N_0)}}{N_0!} \oplus  \frac{A^{(N_0+1)}}{(N_0+1)!} \oplus \cdots  \oplus \frac{A^{(t)}}{t!}\\
           &= \frac{A}{1!}\oplus \frac{A^{(2)}}{2!}\oplus \cdots \oplus \frac{A^{(N_0-1)}}{(N_0-1)!} \\
           & \oplus \frac{A^{(N_0)}}{N_0!} \oplus  \frac{A^{(N_0+1)}}{(N_0+1)!} \oplus \cdots  \oplus \frac{A^{(N_0+p-1)}}{(N_0+p-1)!}\\
           & \oplus \lambda^{(p)}\Big(\frac{A^{(N_0)}}{(N_0+p)!} \oplus  \frac{A^{(N_0+1)}}{(N_0+p+1)!} \oplus \cdots  \oplus \frac{A^{(N_0+p-1)}}{(N_0+2p-1)!}\Big)\\
           & \oplus \lambda^{(2p)}\Big(\frac{A^{(N_0)}}{(N_0+2p)!} \oplus  \frac{A^{(N_0+1)}}{(N_0++2p+1)!} \oplus \cdots  \oplus \frac{A^{(N_0+p-1)}}{(N_0+3p-1)!}\Big)\\
           &\cdots\\
         & \oplus \lambda^{(lp)}\Big(\frac{A^{(N_0)}}{(N_0+lp)!} \oplus  \frac{A^{(N_0+1)}}{(N_0+lp+1)!} \oplus \cdots  \oplus \frac{A^{(N_0+p-1)}}{(N_0+(l+1)p-1)!}\Big)\\
         & \oplus \frac{A^{(N_0+(l+1)p)}}{(N_0+(l+1)p)!} \oplus \frac{A^{(N_0+(l+1)p+1)}}{(N_0+(l+1)p+1)!} \oplus \cdots \oplus  \frac{A^{(t)}}{t!}. 
        \end{aligned}
    \end{equation*}
Where, $l = [(t-N_0+1)/p]-1$, which corresponds to the last complete set of repeating powers, $\Big(\frac{A^{(N_0)}}{k!} \oplus  \frac{A^{(N_0+1)}}{(k+1)!} \oplus \cdots  \oplus \frac{A^{(N_0+p-1)}}{(k+p-1)!}\Big)$ for some $k$. Note that, $l \geq 1$ is guaranteed as $O(A)> 
 N_0+2p$. Now, using simple calculus techniques for finding maxima, we can see that one with the coefficient $\lambda^{(lp)}$ dominates the rest of the terms. Hence, $e^{(A)}$ can be expressed as follows.
\begin{equation*}
        \begin{aligned}
e^{(A)}  &=  \underbrace{ \frac{A}{1!}\oplus \frac{A^{(2)}}{2!} \oplus \cdots \oplus \frac{A^{(N_0-1)}}{(N_0-1)!}}_{R_a}\\ 
& \oplus \lambda^{(lp)}\Big(\frac{A^{(N_0)}}{(N_0+lp)!} \oplus  \frac{A^{(N_0+1)}}{(N_0+lp+1)!} \oplus \cdots  \oplus \frac{A^{(N_0+p-1)}}{(N_0+(l+1)p-1)!}\Big)\\
         & \oplus \frac{A^{(N_0+(l+1)p)}}{(N_0+(l+1)p)!} \oplus \frac{A^{(N_0+(l+1)p+1)}}{(N_0+(l+1)p+1)!} \oplus \cdots \oplus  \frac{A^{(t)}}{t!} 
         \end{aligned}
    \end{equation*}
Also, since $a_{ij}$'s are positive, and $O(A)> N_0 +2p$, in the above increasing series, $R_a$ will be dominated by the rest of the terms, leaving the exponential as follows. 
\begin{equation*}
        \begin{aligned}
e^{(A)}  &=  \lambda^{(lp)}\Big(\frac{A^{(N_0)}}{(N_0+lp)!} \oplus  \frac{A^{(N_0+1)}}{(N_0+lp+1)!} \oplus \cdots  \oplus \frac{A^{(N_0+p-1)}}{(N_0+(l+1)p-1)!}\Big)\\
         & \oplus \frac{A^{(N_0+(l+1)p)}}{(N_0+(l+1)p)!} \oplus \frac{A^{(N_0+(l+1)p+1)}}{(N_0+(l+1)p+1)!} \oplus \cdots \oplus  \frac{A^{(t)}}{t!}\\
        &= \Big( \frac{A^{(N_0+lp)}}{(N_0+lp)!} \oplus \cdots  \oplus \frac{A^{(N_0+(l+1)p)}}{(N_0+(l+1)p)!} \oplus \cdots \oplus  \frac{A^{(t)}}{t!}  \Big)
         \end{aligned}
    \end{equation*}
Thus, all the terms, which are dominating in the power series expansion have power at least $N_0+p$. This gives every column in every dominating term are generalised eigenvectors, and so is their sum.  
    \end{itemize}
\end{proof}
\begin{Example}
    In DEDS models, a Markov chain is often expressed using stochastic matrices, which exhibit recurring relations with respect to the previous state. In this context, it is useful to study the eigenvector and generalised eigenvector problem of matrices such as permutation matrices, stochastic matrices and doubly stochastic matrices. Note that, here we will study about real-stochastic matrices, instead of tropical variants of the same, as they are used in the above-said processes.  \\
    Let $A$ be a real-permutation matrix. Then we have the following:
    \begin{itemize}
        \item [(a)] $\lambda(A)$ = 1.
        \item[(b)] $per(A) = order (\pi)$, where $\pi$ is the corresponding permutation of $A$.
        \item[(c)] Every columns of $e^{(1\otimes A)}=A$ gives generalised eigenvectors of $A$.
    \end{itemize}
    Here, every entry of $A$ is either 0 or 1, and 1 appears exactly once in every column and row. Let $\pi$ be the corresponding permutation in $S_n$, the permutation group of $n$ symbols. Then, each disjoint cycle in $\pi$ denotes critical cycles in $G(A)$.  Then (a) and (b) follow immediately from this structure of $A$.\\
    One can see that, any real-permutation matrix is strictly periodic with the period same as the order of the corresponding permutation. Then, (c) is a consequence of Theorem \ref{gen.eigenvector relation}. 
 Consider the matrix $A = \begin{pmatrix}
     0&1&0&0\\
     1&0&0&0\\
     0&0&0&1\\
     0&0&1&0
 \end{pmatrix}$, with the corresponding permutation is $\pi = (12)(34)$. The disjoint cycles in $\pi$ correspond to different critical cycles. Then, we have, $per(A)= lcm(2,2) = 2$. Here, $e^{(1\otimes A)}= \begin{pmatrix}
     0&1&0&0\\
     1&0&0&0\\
     0&0&0&1\\
     0&0&1&0
 \end{pmatrix} = A$. One can verify that columns of $e^{(1\otimes A)}$ are generalized eigenvectors of $A$ of order 2.    
\end{Example}
\section{declarations}
\subsection{Ethical Approval }
Not Applicable.
\subsection{Availability of supporting data}
Not Applicable.
\subsection{Competing interests}
 The authors have no competing interests as defined by Springer, or other interests that might be perceived to influence the results and/or discussion reported in this paper.
\subsection{Funding}
Not Applicable.
\subsection{Authors' contributions}
H.M formulated the initial research problem to study exponential of Max-plus matrices. A.A and H.M equally contributed towards all results. A.A prepared all figures. A.A and H.M together wrote the manuscript. All authors have reviewed the manuscript.
\subsection{Acknowledgments}
Not Applicable
\bibliographystyle{abbrv}
\bibliography{tropical}

\end{document}